\RequirePackage[english]{babel}
\documentclass{article}
\usepackage{tocloft}

\usepackage[pdftex, linktocpage]{hyperref}
\usepackage{amsmath,amsfonts,amssymb,amsthm,amscd}
\usepackage{comment}

\title{Group covers, $o$-minimality, and categoricity}
\author{A. Berarducci\thanks{Partially supported by PRIN 2007PRYAAF\_004: O-minimalit\`a - Metodi e modelli non standard - Teoria degli insiemi} \\University of Pisa
 \and Y. Peterzil\\University of Haifa \and A. Pillay\thanks{Supported by EPSRC grant EP/F009712/1}\\University of Leeds}
\date{19th September 2010}

\def\R{{\mathbb R}}
\def\NN{{\mathbb N}}
\def\Q{{\mathbb Q}}
\def\Z{{\mathbb Z}}
\def\C{{\mathbb C}}

\def\rf{{{\mathbb R}_{\text{\scriptsize field}}}}
\def\cf{{{\mathbb C}_{\text{\scriptsize field}}}}

\newlength{\tocwidth}
\def\M{M}

\newcommand{\ov}{\overline}

\theoremstyle{plain}
\newtheorem{theorem}{Theorem}
\newtheorem{lemma}[theorem]{Lemma}

\newtheorem{proposition}[theorem]{Proposition}
\newtheorem{corollary}[theorem]{Corollary}

\newtheorem{claim}{Claim}

\newtheorem{question}[theorem]{Question}
\newtheorem{problem}[theorem]{Problem}
\newtheorem{fact}[theorem]{Fact}
\theoremstyle{definition}
\newtheorem{remark}[theorem]{Remark}
\newtheorem{definition}[theorem]{Definition}
\newtheorem{example}[theorem]{Example}
\newtheorem{exercise}[theorem]{Exercise}

\numberwithin{theorem}{section}

\newcommand{\bt}{\begin{theorem}}
\newcommand{\et}{\end{theorem}}

\newcommand{\bl}{\begin{lemma}}
\newcommand{\el}{\end{lemma}}

\newcommand{\bfa}{\begin{fact}}
\newcommand{\efa}{\end{fact}}

\newcommand{\bexa}{\begin{example}}
\newcommand{\eexa}{\end{example}}

\newcommand{\bexe}{\begin{exercise}}
\newcommand{\eexe}{\end{exercise}}

\newcommand{\bprop}{\begin{proposition}}
\newcommand{\eprop}{\end{proposition}}

\newcommand{\bp}{\begin{proof}}
\newcommand{\ep}{\end{proof}}

\newcommand{\bc}{\begin{corollary}}
\newcommand{\ec}{\end{corollary}}

\newcommand{\bd}{\begin{definition}}
\newcommand{\ed}{\end{definition}}

\newcommand{\br}{\begin{remark}}
\newcommand{\er}{\end{remark}}

\newcommand{\bq}{\begin{question}}
\newcommand{\eq}{\end{question}}

\newcommand{\bproblem}{\begin{problem}}
\newcommand{\eproblem}{\end{problem}}

\newcommand{\w}[1]{\widetilde {#1}} 
\newcommand{\ca}[1]{{\mathcal #1}}

\newcommand{\rest}[1]{\hspace {-0.2em} \upharpoonright_{#1}}

\def\i{\iota}
\def\p{\pi}
\def\modZ{{{\rm mod}_\Z}}

\def\M{M}

\newcommand{\ex}[3]{1 \longrightarrow #1 \longrightarrow #2 \longrightarrow #3 \longrightarrow 1}
\newcommand{\exact}[5]{1  \longrightarrow #1 \stackrel{#2}\longrightarrow #3 \stackrel{#4}\longrightarrow #5 \longrightarrow 1}

\newcommand{\esatta}[7]{#1  \longrightarrow #2 \stackrel{#3}\longrightarrow #4 \stackrel{#5}\longrightarrow #6 \longrightarrow #7}
\newcommand{\structure}[5]{(#1,#2,#3, #4, #5)}
\newcommand{\str}[3]{(#1,#2,#3)}

\begin{document} 

\maketitle

\begin{abstract} We study the model theory of covers of groups definable in $o$-minimal structures. This includes the case of covers of compact real Lie groups. In particular we study categoricity questions, pointing out some notable differences with the case of covers of complex algebraic groups studied by Zilber and his students. We also discuss from a model theoretic point of view the following question, related to ``Milnor's conjecture'' in \cite{Milnor:83}: is every finite central extension of a compact Lie group isomorphic to a topological extension? 
\end{abstract} 

\begin{center}
\begin{minipage}{\tocwidth}
\small
\tableofcontents
\end{minipage}
\end{center}

\section{Introduction and preliminaries} 
There are at least two inspirations for the current paper. The first is the work of Zilber and students, \cite{Zi:06},
\cite{Gav:06}, \cite{Bays:09}, around the model theory and categoricity, sometimes infinitary,  of universal covers of commutative complex algebraic groups. We were interested in studying analogous questions for real Lie groups or groups definable in $o$-minimal structures.  The second is the paper \cite{HrPePi:08b} by the second and third authors together with Hrushovski, which studied group extensions definable in an $o$-minimal structure, as well as topological covers of definable Lie groups. In section 2 of this paper we study abstract finite central extensions  $\quad \exact \Gamma \iota G \pi H$  where $H$ is, say, definable in an 
$o$-minimal expansion $M$ of the field of real numbers, but now no ``tameness"  assumption is made on $G$. We formulated the following conjecture, which we were subsequently informed to be part of ``Milnor's conjecture"  \cite{Milnor:83}:
\newline
(*) $G$ can be equipped with Lie group structure making $\pi: G\to H$ a topological covering.
\newline
As far as we understand the  case  when $H$ is compact (so a compact Lie group) is still open. We give, in Proposition 2.2 and Theorem 2.3, some model-theoretic equivalences to (*). We give suitable formulations when $M$ is just an $o$-minimal expansion of an arbitrary real closed field, as well as a positive solution in this context, for commutative $H$. We also make some comments on the (definably) compact and semisimple cases. Among the model-theoretic equivalences to (*) that we mention, is the {\em stable embeddability} of the structure $M$ in the structure $((G,\cdot),\pi,M)$ . In the case where $H$ is compact, we highlight the role of the ``infinitesimal subgroup" $H'^{00}$ of a saturated elementary extension $H'$ of $H$. There is some overlap in our section 2 with current work by Edmundo, Jones, and Peatfield  \cite{EdJoPe:10} and we will give precise references. 

Section 3 of the paper deals with categoricity issues and universal covers $\pi:G \to H$ of definable Lie groups $H$. We make use of results from $\cite{HrPePi:08b}$ to prove the following strong relative categoricity statement (Theorem 3.4): there is a single $L_{\omega_{1},\omega}$-sentence $\sigma$ true of $((G,\cdot),\pi, M)$  such that if $(G_{1},M_{1})$, $(G_{2},M_{2})$ are models of $\sigma$ then any isomorphism between the real closed fields $M_{1}$ and $M_{2}$ lifts to an isomorphism between the two structures. On the other hand it does not suffice that $\sigma$ simply fixes the isomorphism type of  $\ker(\pi)$.  We go through the case of $H = \R/\Z$ in detail, summarizing the situation in Theorem 3.12. 
We also discuss stable embeddability, and ask questions (such as 3.14) around suitably generalizing the material from section 2 to the context of arbitrary finitely generated central extensions of definable Lie groups. 

\vspace{5mm}
\noindent
Our model theoretic notation is quite standard. We normally feel free to work in $M^{eq}$  (or $T^{eq}$) so by a definable set in a structure $M$ we mean a set definable (possibly with parameters) in $M^{eq}$. On the other hand in cases such as $o$-minimal structures, there is a privileged sort. For certain technical reasons related to equipping a definable group with a ``manifold topology", when we speak of a group $G$ definable in an $o$-minimal structure $M$, we will tacitly assume that in fact the universe of $G$ is a subset of some $M^{n}$. Of course when the (one-sorted) $o$-minimal structure is an expansion of an ordered group then it eliminates imaginaries, so up to definable bijection this is no restriction. We often will work in a ``saturated model" which means a sufficiently saturated and homogeneous model. For a set $A$ of parameters from a structure $M$, by an $A$-definable set we mean a set definable in $M$ with parameters from $A$.

In Section 2 we will use the notion of a definable set in an ambient structure being $o$-minimal. The meaning is as follows: Let $N$ be a structure, 
and $X$ a $\emptyset$-definable set in $N$ equipped with a linear ordering $<$ also $\emptyset$-definable in $N$. We say that $X$  (or $(X,<)$) is {\em $o$-minimal} in $N$ if every subset of $X$ which is definable in $N$ is a finite union of intervals (with endpoints in $X\cup \{+\infty, -\infty\}$) and points. We say that $X$ is {\em strongly $o$-minimal} in $N$, if the same is true in a saturated elementary extension of $N$ (equivalently the number of intervals and points is bounded inside definable, in $N$, families of definable subsets of $X$).

Notions such as {\em interpretability}, and {\em stable embeddability}, figure prominently in the paper, so we take the opportunity here to  fix our understanding of these notions, at least with respect to the key structures that concern us.

\vspace{2mm}
\noindent
Consider an exact sequence of groups

$$E : =  \quad \exact \Gamma \iota G \pi H $$

We may consider $E$ as a 3-sorted structure $\structure \Gamma \i G \p H$ with sorts for $\Gamma,G,H$ and functions for the group operations and the homomorphisms $\iota, \pi$. Let us call this the {\em pure group language}. Now suppose that the group $H$ is definable, let's say without parameters, in some first order structure $\M$ (which could be $H$ itself). We may then consider the richer structure $N := \structure \Gamma \i G \p \M$ where we have omitted the $H$-sort and added a sort for $\M$ (from these data we can recover the group $H= \pi(G) \cong G/\i(\Gamma)$). Since $\Gamma$ is isomorphic to the kernel of $\pi$, we will often use the abridged notation $$N = \str G \p M,$$ but this is only a matter of notational convenience since in practice it is convenient to allow $\i$ not to be the inclusion. So officially $\str G \p M$ has sorts for $\Gamma, G, \M$ and symbols for the group operation of $G$, the homomorphisms $\i$ and $\p$, and the relations and functions of $\M$. When $\M$ is just the group $H$ we are in the pure group language. An isomorphism between two such structures $N = \str G \p \M$ and $N' = \str {G'} {\p'} {\M'}$ is given by a commutative diagram with exact rows

\[
\begin{CD}
1 @>>> {\Gamma} @>>> {G} @>>> {H} @>>> 1 \\
@. @VV{f_{\Gamma}}V @VV{f_{G}}V @VV{f_{H}}V  \\
1 @>>> {\Gamma'} @>>> {G'} @>>> {H'} @>>> 1 
\end{CD}
\]
where $f_{\Gamma},f_G,f_{H}$ are group isomorphisms and $f_H$ is induced by an isomorphism from $\M$ to $\M'$.

\vspace{2mm}
\noindent
In this situation, with $N = (G,\pi,M)$, we will say that {\em $N$ is naturally interpretable in $M$} (without parameters) if there are, in the structure 
$M$,  a $\emptyset$-definable group  $G_{0}$  and a $\emptyset$-definable homomorphism  $\pi_{0}:G_{0}\to H$, 
such that the structure $N = \str G \p \M$ is isomorphic to the structure $N_{0} = (G_{0},\pi_{0},M)$ via an isomorphism which is the identity on $M$. If $G_{0}$ and  $\pi_{0}$ are allowed to be definable with parameters in $M$ we will say that $N$ is naturally interpretable in $M$ with parameters, maybe mentioning the parameters explicitly.
So the main point about naturality here is that the interpretation of $N$ in $M$ should be ``over $M$". 

We will also speak of  $N := \structure \Gamma \i G \p \M$  being (naturally) interpreted (without parameters, with parameters) in the $2$-sorted structure $(\Gamma, M)$, (where $\Gamma$ is equipped just with its group structure), and again this refers to an interpretation {\em over}  $\Gamma$ and $M$. 

In our context we will say that $M$ is {\em stably embedded} in $N$ if any subset $X$ of $M^{n}$ which is definable with parameters in the structure $N$, is definable with parameters in the structure $M$. A related notion (which sometimes, although not in this paper, is taken as part of the definition of stably embedded) is that any subset $X$ of $M^{n}$ which is definable, without parameters, in $N$ is definable, without parameters, in $M$. The latter property  passes to any model $N' = (G',\pi',M')$ elementarily equivalent to $N$. But  stable embeddability does {\em not} have to be preserved by elementary equivalence. 

It is worth noting that if  $N$ is (naturally) interpretable in $M$ then such an interpretation is reflected in $Th(M)$ in the sense that $Th(M)$ includes the sentences expressing that  $(G_{0},\pi_{0},M)$ is a model of $Th(N)$. So for any model $M'$ of $Th(M)$, $(G_{0}(M'),\pi_{0}(M'),M')$ is a model of $Th(N)$ (where $G_{0}(M')$ etc... denote the interpretations of the relevant formulas in $M'$). But it is {\em not}  necessarily the case that for any model $N' = (G',\pi',M')$ of $Th(N)$, $N'$ is isomorphic to $(G_{0}(M'),\pi_{0}(M'),M')$.

\vspace{5mm}
\noindent
 Typically $\M$ will be the complex field $\cf$, or the real field $\rf$. 
If $H$ is definable in $\cf$ it is also definable in $\rf$ (identifying $\C$ with $\R\times \R$), but the language of $\str G \p \rf$ is richer than the language of $\str G \p \cf$, so the model-theoretic properties of these two structures are very different. 
In particular the complex field $\cf$ has many automorphisms, while the real field $\rf$ is rigid. The following 
example is instructive.

\bexa  Consider the exact sequence $$E:= \quad 0\longrightarrow \Z\stackrel \i \longrightarrow (\C,+) \stackrel{\exp}\longrightarrow {\C^*} \to 1$$ given by the complex exponential function, viewed as above as a structure $(V,\pi,\cf)$, say. Let $(G,+)$ be an abelian divisible torsion free abelian group. Zilber 
\cite{Zi:06} proved that any exact sequence of the form $$\quad 0\longrightarrow \Z\stackrel \i \longrightarrow (G,+) \stackrel{\exp}\longrightarrow \C^* \to 1$$ is isomorphic to $E$ via an isomorphism which can be chosen to be the identity on the $\Gamma$-sort ($=\Z$). However one cannot require that the isomorphism is also the identity on the $H$ (i.e. field) sort. So the relative categoricity statement fails when we replace the complex field $\cf$ by the real field $\rf$. 

Another difference between the complex and real setting is as follows (and was already mentioned in \cite{HrPePi:08b}):
The structure  $(V,\pi,\cf)$ above is superstable of finite $U$-rank. It is neither  naturally nor unnaturally interpretable in the $2$-sorted  superstable structure $(\Z, \cf)$, which we explain now. The $2$ sorts  $(\Z,+)$  and $\cf$  are orthogonal definable groups of $U$-rank $1$. Stability theoretic arguments yield that any group $G$ definable in  $(\Z, \cf)$, is definably an almost direct product of a group $G_{1}$ definable in $(\Z,+)$ and a group $G_{2}$ definable in $\cf$. So if $(V,\pi,\cf)$ were interpretable in $(\Z, \cf)$ then the associated exact sequence would definably split. But it does not even abstractly split. 
On the other hand, from \cite{HrPePi:08b} it follows that $(V,\pi,\cf)$ {\em is}  ``naturally" interpretable in 
the structure $(\Z,\rf)$. 
\eexa

\vspace{2mm}
\noindent
All three authors would like to thank the Durham Symposium on New Directions in the Model Theory of Fields (July 2009) where discussions on the topic of this paper began.
The first and third authors would like to thank the G\"odel Centre at the University of Vienna for its hospitality in June 2010 when they were invited professors and continued work on the paper.

\section{Finite central extensions}

Let $\M$ be an o-minimal structure, let $H$ be a definably connected group in $\M$ and let $\pi \colon G \to H$ be a finite central extension of $H$. We study the following problem:

\bproblem 
\label{problem} 
When is  $(G,\pi, M)$ (naturally) interpretable in $\M$?
\eproblem

It turns out that in this context ``interpretable'' is equivalent to ``naturally interpretable'', but for the proof of Theorem \ref{ominimal} (see below) it will be convenient to work with natural interpretations. The equivalence of the two notions follows immediately from the fact that, by results in \cite{OtPePi:96}, in an $o$-minimal expansion of a real closed field, any definable real closed field is definably isomorphic to the ground field.
 
Let us begin with an observation implicit in the literature. 

\bprop \label{real-case} Assume $\M$ is an o-minimal expansion of the reals. Then $(G,\pi,M)$ is interpretable in $\M$ if and only if  $G$ can be given a group topology that makes $\pi$ a topological covering (hence a homomorphism of Lie groups).
\eprop

\bp For the left implies right direction: Via \cite{Pi:88}  $G$ is definably equipped with the structure of a real Lie group  and $\pi:G\to H$ is a (definable) homomorphism of Lie groups, and as the kernel is finite, must be a topological covering. The right implies left direction appears in \cite{HrPePi:08b} (see Theorems 2.8 and 8.4) as well as in \cite{EdJoPe:10}  (Theorem 1.4). 
\ep

So to address our question we need to compare abstract group extensions with {\em topological group extensions}, meaning the extension is a topological covering (which is stronger than simply being a quotient map). Note that a covering homomorphism of connected Lie groups has central kernel, so the centrality of the kernel is a necessary condition to have a topological extension. The question whether any abstract finite central extension $\ex \Gamma G H$ of a compact connected real Lie group $H$ is equivalent to a topological extension (i.e. to an extension of Lie groups) is related to the ``Friedlander-Milnor conjecture'' or just ``Milnor conjecture" \cite{Milnor:83}. We do not know the answer even for $H = SO_3(\R)$ and $\Gamma = \Z/2\Z$.

\subsection{Stable embeddedness}
Let $\M$ be an o-minimal structure. Let $H$ be a definably connected group in $\M$ and let $\esatta 1 \Gamma \i G \pi H 1$ be a finite central extension. We will assume, here and subsequently, that $H$ is $\emptyset$-definable in $M$ (by adjoining constants if need be). 

\bt \label{ominimal} The following are equivalent. 
\begin{enumerate}
\item $\str G \p \M$ is naturally interpretable in $\M$; 
\item $\M$ is stably embedded in $\str G \p \M$;
\item $\M$ is an  $o$-minimal set in $\str G \p \M$. 
\end{enumerate}
Moreover these conditions imply that  $M$ is strongly $o$-minimal in $\str G \p \M$ and that the stable embeddability statement holds in all models of $Th(\str G \p \M)$.
\et

We do not know whether {\em any} finite central extension of $H$ satisfies the above conditions, but we will prove later that this holds in the abelian case. 

Again much of the proof of Theorem 2.3 is contained in earlier papers. For example $2\implies 1$ can be seen to follow from Theorem 8.2 of \cite{HrPePi:08b}. And $3\implies 1$ can be seen as a restatement of Corollary 1.2 of \cite{EdJoPe:10}. Neverthless for the benefit of the reader we will give more or less direct proofs, following a sequence of lemmas.

\bl \label{section} 
\begin{enumerate}
\item There is a section $s \colon H \to G$ of $\pi$ which is $\emptyset$-definable in $(G,\M,\pi)$. 
\item The extension $\pi\colon G \to H$ is definably isomorphic in $\str G \p M$ to an extension $\pi'\colon G' \to H$ such that the underlying set of $G'$ and the homomorphism $\pi'$ are definable in $M$, and the group operation of $G'$ is definable in $\str G \p M$. 
\end{enumerate}
\el

\bp (1) Let $n = |\ker (\pi)|$. Since $H$ is definably connected, by \cite{HrPePi:08b} (Lemma 8.1(ii)) there is some $k$ such that every $y\in H$ can be written in the form $y_1^n \cdot \ldots \cdot y_k^n$. By definable choice there are definable functions $r_1,\ldots, r_k \colon H \to H$ such that for all $y\in H$ we have $y = r_1(y)^n \cdot \ldots \cdot r_k(y)^n$. Since the extension is central and $|\ker(\pi)|=n$, any two elements $u,v\in G$ with $\pi(u)=\pi(v)$ satisfy $u^n=v^n$. So we can define a section $s\colon H\to G$ as follows. Given $y\in H$ pick $x_1\in \pi^{-1}(r_1(y)), \ldots, x_k \in \pi^{-1}(r_k(y))$ and define $s(y) = x_1^n \cdot \ldots \cdot x_k^n$. This proves the first part.

(2) We can use $s$ to define a bijection $f\colon H \times \Gamma \to G$ sending $(x,a)$ to $s(x)a$. After fixing finitely many constants from $\M$, we can assume that $\Gamma \subset \M$. So the set $H\times \Gamma$ is definable in $\M$. Note that $\pi':= \pi \circ f \colon H\times \Gamma \to H$ is the projection on the first coordinate. Put on $H\times \Gamma$ the unique group structure making $f$ into an isomorphism, and call $G'$ the resulting group. 
\ep 

\bc \label{tame-int} Let $\pi \colon G \to H$ be a finite central extension of a definably connected group $H$. Assume that $\M$ is stably embedded in $\str G \p \M$. Then $\pi \colon G \to H$ is naturally interpretable in $\M$. 
\ec

\bp By Lemma \ref{section} we can assume that both the domain of $G$ and the map $\pi$ are definable in $\M$, while the group operation $\mu$ on $G$ is definable in $\str G \p \M$. Since $\M$ is assumed to be stably embedded in $\str G \p \M$, $\mu$ must be definable in $\M$. 
\ep  

The assumption that $\M$ is stably embedded can be weakened: it suffices that $\M$ be $o$-minimal  in $\str G \p \M$. We need:

\bl \label{good-reduction} Let $\M$ be an $o$-minimal structure in a language $L$, and let $L_{0}$ be a sublanguage of $L$ including the ordering on $M$. Let $\pi\colon G \to H$ be a group homomorphism with finite kernel $\Gamma$, all definable in $\M$. Assume that $dom(G)$ (the underlying set of $G$) as well as the map $\pi$ and the group $H$ are all definable in the language $L_{0}$.  Then also the group operation of $G$ is definable in the language $L_{0}$, whence the group $G$ is definable in the reduct $\M|L_{0}$ of $\M$. \el

\bp 
The idea is to show that $G$ has an $L_0$-definable topology, and then to observe that the group structure is determined by the topology. To this aim fix a finite cover $\ca F$ of $H$ by $L_0$-definable definably simply connected open sets. For $U\in \ca F$, each definably connected component of $\pi^{-1}(U)$ is $L_0$-definable and homeomorphic to $U$ via the projection $\pi$. We have thus endowed $G$ with the structure of a definable manifold with an $L_0$-definable atlas, making $\pi \colon G\to H$ into an $L_0$-definable covering map. Now the $L_0$-definable group operation of $H$ can be lifted uniquely (via uniform lifting of paths) to an $L_0$-definable group operation $\mu$ on the $L_0$-definable covering space $G$, making $\pi$ into a homomorphism. This $\mu$ must then coincide with the original group operation on $G$. \ep

\bc \label{interpretable} Assume that $\M$ is $o$-minimal in $\str G \p \M$. Then $\str G \p \M$ is interpretable in $\M$. \ec

\bp By Proposition \ref{section} we can assume that both the domain of $G$ and the map $\pi$ are definable in $\M$. Let $\M_{ind}$ be the expansion of $\M$ obtained by adding symbols for all relations between elements of $\M$ which are $\emptyset$-definable in $\str G \p \M$. So the group operation $\mu$ of $G$ is $\emptyset$-definable in $\M_{ind}$. By our assumption, $\M_{ind}$ is o-minimal and 
$\M$ is the reduct of $\M_{ind}$. So by Lemma \ref{good-reduction} $\mu$ is definable in the structure $\M$. 
\ep 

Theorem \ref{ominimal} follows from the above: Firstly $2 \rightarrow 1$ is Corollary \ref{tame-int}, $3 \rightarrow 1$ is Corollary \ref{interpretable}, and clearly $1$ implies each of $2$ and $3$ (using the naturality of the interpretation). The moreover statement is also clear: Indeed from Lemmas \ref{tame-int}, \ref{good-reduction}, \ref{interpretable} and their proofs, it follows that the equivalent conditions of Theorem \ref{ominimal} imply the following strengthening of point 1: \\
(1*) $(G,\pi,M)$ is naturally interpretable in $M$ ``definably in $(G,\pi,M)$''
 (namely the isomorphism of the interpretation is definable).\\
To deduce from (1*) the ``moreover part'' of \ref{ominimal} one can simply observe that 
natural {\em definable} interpretations pass to all models of the theory of $(G,\pi,M)$, 
so we get the corresponding strong forms of 2 and 3.

\subsection{The abelian case}
As in the previous subsection
let $H$ be a definably connected group $\emptyset$-definable in the $o$-minimal structure $\M$ and let $\esatta 1 \Gamma \i G \pi H 1$ be a finite central extension. In this subsection we make the additional assumption that $H$ is abelian, and give a positive solution to Problem \ref{problem}. 

\bl \label{bilinear} $G$ is abelian.
\el
\begin{proof}
So $\pi:G\to H$ is an extension of $H$ by $\Gamma < Z(G)$. Note that for $x,y\in G$ and $c\in Z(G)$ we have $[xc,yc]=[x,y]$. Since the extension $\pi\colon G \to H$ is central, we then have a well defined map $b: H\times H \to \Gamma$, given by $b(\pi(x),\pi(y)) = [x,y]$. For $a\in H$, $\{t\in H:b(a,t) = 0\}$ is a subgroup of $H$ (since in any group $[x,yz]=[x,y][x,z]^y$) which has finite index if $\Gamma$ is finite. But $H$, as a definably connected commutative group in an $o$-minimal structure, is known to be divisible, so has no proper subgroups of finite index. Thus $G$ is commutative.
\end{proof}

\bt \label{abelian} $(G, \pi, \M)$ is naturally interpretable in $\M$.
\et

\bp 
Let $|\Gamma| = n$.
By Lemma \ref{bilinear} $G$ is abelian. Assume first that  $nG = G$. 
There a surjective group homomorphism $\varphi \colon H \to nG$ sending $\pi(g)$ to $ng$. 
Its kernel is $L = \pi(G[n]) \subset H[n]$, a finite subgroup of $H$ (since any definable abelian group has finitely many elements of any given order). 
Let $\varphi_L \colon H/L \to nG$ be the induced isomorphism. Composing with $\pi\colon G \to H$ we get a surjective group
homomorphism $\pi'\colon H/L \to H$, $\pi'(x+L) = nx$, which is definable (without parameters) in $M$ and isomorphic over $H$ to $\pi\colon G \to H$. 

Now consider the case where $nG \neq G$. There is a surjective homomorphism $nG \times \Gamma \to G$, $(x,c)\mapsto x + c$, with finite kernel. So $G$ is the almost direct product of $\Gamma$ and $nG$. Exactly as in the 
previous paragraph, $(nG,\pi|nG, M)$ is naturally interpretable in $M$. As $\Gamma$ is finite, also $\pi \colon G\to H$ is interpretable in $M$.  
\ep

\subsection{The ``infinitesimal" subgroup}

We now specialize Problem \ref{problem} to the case when $H$ is definably compact. In this case the 
``intrinsic infinitesimal subgroup'' $H^{00}$ of $H$ (see \cite{Pi:04,BeOtPePi:05,HrPePi:08}) gives us some valuable information. So for this subsection, let $\M$ be a saturated o-minimal expansion of a real closed field, let $H$ be a definably compact definably connected definable group in $\M$ and let $$\esatta 1 \Gamma \i G \pi H 1$$ be a finite central extension of $H$ with $|\Gamma|=n$. We do not always use the fact that $\M$ expands a real closed field, but it is convenient. We will prove, in this section the following:

\bt \label{split} The following are equivalent. 
\begin{enumerate}
\item $\str G \p \M$ is naturally interpretable in $\M$.
\item The set $G_{00}$ consisting of $nth$-powers of elements of $\pi^{-1}(H^{00})$ is a subgroup of $G$. 
\item The sequence $\esatta 1 \Gamma \i G \p H 1$ splits over the infinitesimal subgroup $H^{00}$ of $H$ (namely $\esatta 1 \Gamma \i {\pi^{-1}(H^{00})} \p {H^{00}} 1$ splits as an extension of abstract groups). 
\end{enumerate} 
\et 

\br If $\M$ is not saturated (for instance $\M$ is the real field), we can still apply the theorem by passing to a saturated extension. So (2) holds in a saturated extension if and only if (1) holds in the original model (using the ``moreover part'' of Theorem \ref{ominimal}). \er

\bl \label{unique-div} $H^{00}$ is uniquely divisible in the following sense: for every $x\in H^{00}$ and every $n\in \NN^*$, there is unique $y\in H^{00}$ such that $y^n=x$. 
\el
\bp 
We already know that $H^{00}$ is divisible and torsion free (by Theorem 4.6 of \cite{Be:07}). In the abelian case this implies unique divisibility. In the general case, by \cite{HrPePi:08b} we can write $H$ as an almost direct product $H = A\times_\Gamma B$ of an abelian definably connected group $A$ and a semisimple definably connected group $B$ (so the intersection $\Gamma = A\cap B$ is finite and central, and each element of $A$ commutes with each element of $B$). We have $H^{00} = A^{00}B^{00}$. Since $H^{00}$ is torsion free and every element of $A\cap B = \Gamma$ is torsion, the intersection $A^{00}\cap B^{00}$ reduces to the identity element, namely $H^{00}$ is the direct product $A^{00} \times B^{00}$. So it only remains to consider the case when $H$ is semisimple. In this case, replacing $H$ with a definably isomorphic group, we can assume that $H$ is definable without parameters in the pure field language (again by results in \cite{HrPePi:08b}). So it makes sense to consider $H(\R)$, and by \cite{Pi:04} we have $H/H^{00} \cong H(\R)$ with the natural projection $H(\M)\to H/H^{00}$ corresponding to the standard part map $st: H(\M)\to H(\R)$. The exponential map $\exp \colon T_e(H(\R)) \to H(\R)$ is a local homeomorphism, so we can fix a small convex neighborhood $U$ of $T_e(H(\R))$ such that $\exp|U$ is a homeomorphism onto its image $V\subset H(\R)$. This $V$ is a ``uniquely divisible'' open neighborhood of the identity, in the sense that for each $x\in V$ and each $n$ there is a unique $y\in V$ with $y^n=x$ and $y^{i}\in V$ for all $i\leq n$. Unfortunately $V$ is not definable in the field language, since we used $\exp$, however it can be approximated by an open set $V'$ which is definable without parameters in the pure field language and it is still uniquely divisible (take any semialgebraic set $V'$ between $V$ and $V^{1/2}$, where $V^{1/2} = \exp(\{x\in U : 2x \in U\})$). By completeness of the theory of real closed fields, $V'(\M)\subset H(\M)$ remains uniquely divisible. But $H^{00} \subset V'(\M)$, so $H^{00}$ is uniquely divisible. 
\ep

\bc \label{powers}
Let $G_{00}$ be the set of $nth$-powers of elements of $\pi^{-1}(H^{00})$ (where  $n = |\Gamma|$). 
Then $\pi \rest {G_{00}} \colon G_{00} \to H^{00}$ is bijective. 
\ec
\bp It is surjective because $H^{00}$ is divisible, and injective because $H^{00}$ is uniquely divisible. \ep

\bp[Proof of Theorem \ref{split}] Thanks to Corollary \ref{powers} (and the fact that a group homomorphism maps $nth$-powers to $nth$-powers), the only possible splitting homomorphism $s\colon H^{00}\to \pi^{-1}(H^{00})$ of $\esatta 1 \Gamma \i {\pi^{-1}(H^{00})} \p {H^{00}} 1$ is given by the inverse of $\pi\rest{G_{00}}$ (which is a group homomorphism if and only if $G_{00}$ is a group). This proves that (2) is equivalent to (3). 

(1$\to$ 2) Suppose $\pi \colon G \to H$ is (naturally) interpretable in $\M$. We must show that $G_{00}$ is a group. 
So we assume that $G,\pi$ are definable in $\M$. There is no harm in assuming that $G$ is definably connected. As $H$ is definably compact and $\ker(\pi)$ is finite, $G$ is also definably compact. We will show that  the ``intrinsic infinitesimal subgroup" $G^{00}$ of $G$ coincides with what we have called above $G_{00}$, which suffices to give (2). We know that $G^{00}$ is ``torsion-free" so has trivial intersection with $\Gamma$. It also clearly maps surjectively to $H^{00}$ under $\pi$ and hence $\pi|G^{00}$ is an isomorphism between $G^{00}$ and $H^{00}$.  
As $G^{00}$ is divisible $G^{00}$ coincides with the set of $nth$ powers of its elements, and hence $G^{00}\subseteq G_{00}$. By Corollary 2.13 we see that $G^{00} = G_{00}$ as required. 

(2$\to$ 1) We want to prove the (natural) interpretability of $(G,\pi,M)$ in $M$ (assuming that $G_{00}$ is a group). We first make a convenient reduction:
\newline
{\em Claim.} We may assume that $G$ is definably connected in the structure  $(G,\pi,M)$, namely has no proper definable subgroup of finite index.
\bp  First note that the definably connected definably compact group $H$ has no proper subgroups of finite index
(as for example each element of $H$ is an $nth$ power for all $n$). So if $G_{1}$ is a subgroup of $G$ of finite index (without loss normal) then $\pi(G_{1}) = H$, and the index of $G_{1}$ in $G$ is bounded by $|\Gamma|$. It follows that there is a smallest definable (in $(G,\pi,M)$) subgroup of $G$ of finite index, $G^{0}$ say. $G^{0}$ maps onto $H$ under $\pi$ and  $G$ is a quotient of $G_{0}\times \Gamma$ by a finite subgroup. So if $(G^{0},\pi|G^{0},M)$ is (naturally) interpretable in $M$ so is $(G,\pi,M)$.
\ep

We continue with the proof of 2$\to$ 1. Assume that $G_{00}$ is a group. Then the restriction of $\pi$ to $G_{00}$ is a group homomorphism and by Corollary \ref{powers} it is an isomorphism onto $H^{00}$ (so $G_{00}$ is divisible and torsion free). We must show that $\str G \p \M$ is interpretable in $\M$. Note that $G_{00}$ is type-definable in the structure $\str G \p \M$ and has bounded index in $G$ because $\pi$ induces a morphism $\pi_1 \colon G/G_{00} \to H/H^{00}$ with a finite kernel (isomorphic to $\Gamma$).  So if we put on $G/G_{00}$ the logic topology, $\pi_1$ is a continuous homomorphism of connected compact groups. (Connectedness of $G/G_{00}$ is because of definable connectedness of $G$.)  Since its kernel is finite and $H/H^{00}$ is a Lie group, we easily conclude (see below) than $\pi_1$ is a local homeomorphism (and therefore $G/G_{00}$ is also a Lie group). Granted this claim, we continue as follows. By \cite[Prop. 8.3]{BeMa:10} every connected finite extension of $H/H^{00}$ in the category of Lie groups, comes from a definable extension of $H$. Applying this result to the extension $\pi_1 \colon G/G_{00}\to H/H^{00}$, we obtain a definable group extension $\pi' \colon L \to H$ and an isomorphism $f \colon L/L^{00} \cong G/G_{00}$. To prove that the abstract extension $\pi\colon G \to H$ is interpretable in $\M$ it suffices to show that it is isomorphic to the definable extension $\pi'\colon L \to H$. To this aim define $\phi \colon L \to G$ as the map which sends $x\in L$ to the unique $y\in G$ such that $f(xL^{00}) = yG_{00}$ and $\pi(x) = \pi'(y)$. It is easy to see that $\phi$ is indeed an isomorphism of group extensions. It remains to prove the missing claim above, namely that the continuous homomorphism $\pi_1\colon G/G_{00}\to H/H^{00}$ is a local homeomorphism. Clearly it is a closed map since $G/G^{00}$ is compact. It is also an open map being a continuous homomorphism of topological groups. Now, since $\ker(\pi_1)$ is finite, for any sufficiently small open neighbourhood $O \subset G/G_{00}$ of the identity we have $OO^{-1} \cap \ker(\pi_1) = \emptyset$. So $\pi_1$ is locally injective.
\ep

\subsection{The semi-simple case} 

Let us consider again Problem \ref{problem} when $H$ is definably compact and definably connected. By Theorem \ref{split} $\str G \p \M$ is interpretable in $\M$ if and only if the sequence $$\esatta 1 \Gamma \i G \p H 1$$ splits over $H^{00}$. By \cite{HrPePi:08b} $[H,H]$ is a definably semisimple definable group, $H = Z(H)[H,H]$, and $Z(H)\cap [H,H]$ is finite. It follows that $H^{00}$ is the direct product of $Z(H)^{00}$ and $[H,H]^{00}$. So the sequence splits over $H^{00}$ if and only if it splits over $Z(H)^{00}$ and over $[H,H]^{00}$. The former condition is always true by the results of subsection 2.2, (Theorem \ref{abelian}). So we have reduced our problem to the case when $H$ is semisimple. Note that in this case $H$ is perfect (i.e. $H=[H,H]$). Unfortunately we are not able to carry out a complete analysis of this case, but we have the following partial result. 

\bprop \label{commutators} Let $H$ be definably compact definably connected and semisimple. Let $K= \pi^{-1}(H^{00})$. 
\begin{enumerate}
\item If $[K,K] \cap \Gamma = 1$, then the extension $\pi\colon G \to H$ splits over $H^{00}$ (so $\pi \colon G \to H$ is interpretable in $\M$). 
\item In any case we have $[K,K]_1 \cap \Gamma = 1$, where $[K,K]_1$ is the set of commutators of $ K$. 
\item So $[K,K] = [K,K]_1$ is a sufficient condition for the interpretability of $\pi\colon G \to H$ in $\M$. 
\end{enumerate}
\eprop

The proof is given later. 
Note that if $H$ is definably compact, definably connected, and semisimple, then $H=[H,H]=[H,H]_1$ (for real Lie groups this is Goto's theorem). We can also show that $H^{00} = [H^{00},H^{00}]$, but we do not know whether this is equal to $[H^{00},H^{00}]_1$. (This would be the case if and only if the following were true: the commutator map $[,]\colon H\times H \to H$ sends each neighbourhood of the identity in $H\times H$ to a neighbourhood of the identity in $H$.) In any case what we need is stronger, namely $[K,K]=[K,K]_1$. 

\bl \label{centralizer} Given a definably compact group $H$, if $A$ is an abelian subgroup of $H^{00}$, then there is a definably connected definable abelian subgroup $L$ of $H$ which contains $A$. \el 
\bp Note that $A\subset Z(C_{H}(A))$ (the centre of the centralizer in $H$ of $A$) and that $Z(C_H(A))$ is definable (even if $A$ may not be). Note that $Z(C_{H}(A))^{00}$ is contained in $Z(C_{H}(A)) \cap H^{00}$. But the latter is divisible and torsion-free, hence by Corollary 1.2 of \cite{BeOtPePi:05} for example, $Z(C_{H}(A))^{00} = Z(C_{H}(A)) \cap H^{00}$. As  $Z(C_{H}(A))^{00} \leq Z(C_{H}(A))^{0} \leq Z(C_{H}(A))$, we see that $A$ is contained in $Z(C_{H}(A))^{0}$, as required.
\ep 

\bp[Proof of Proposition \ref{commutators}]
(1) Clearly $[K,K]$ projects onto $[H^{00},H^{00}] = H^{00}$. If $\Gamma \cap [K,K] = 1$ then $\pi\rest{[K,K]}$ is injective, so the extension $\pi \colon G \to H$ splits over $H^{00}$. 

(2) Let $\gamma \in \Gamma \cap [K,K]_1$. So we can write $\gamma = [a,b]$ for some $a,b\in K$. We must prove that $\gamma = 1$. Since $\pi([a,b]) = \pi(\gamma) = 1$, it follows that $x:= \pi(a)$ commutes with $y:= \pi(b)$. By Lemma \ref{centralizer} there is an abelian definably connected definable subgroup $L$ of $H$ containing $x$ and $y$. By Lemma \ref{bilinear} $\pi^{-1}(L)$ is abelian. So $[a,b]=\gamma = 1$. 
\ep 

\section{Universal covers}

\subsection{2-cocycles and sections}

Let $H$ be a definable group in a structure $\M$. 
Let $$\esatta 1 \Gamma \i G \p H 1$$ be a central extension of $H$, namely an exact sequence of groups with $\Gamma \leq Z(G)$ (identifying $\Gamma$ with $\ker(\pi)$). Note that we do not assume $\Gamma$ to be finite, so the results of the previous sections do not apply.
Let $s\colon H \to G$ be a {\em section} of $\pi$, namely $s$ a function such that $\pi\circ s$ is the identity on $H$. There is a bijection $f\colon H\times \Gamma\to G$ sending $(x,c)$ to $s(x)\cdot c$, and we can put on $H\times \Gamma$ the unique group operation making $f$ into an isomorphism. This group operation can be described explicitly as follows. Consider the  
{\em 2-cocycle} $h\colon H \times H \to \Gamma, (x,y) \mapsto s(xy)^{-1}s(x)s(y)$ induced by $s$. Then on $H\times \Gamma$ we have the following group operation: 
\begin{equation}\label{inter} (x,c)(y,d)=(xy,c+d+h(x,y)),
\end{equation}
where we have written the group operation on $\Gamma$ additively. Call $H\times_h \Gamma$ the resulting group. Then $H\times_h \Gamma \cong G$ via $f\colon (x,c)\mapsto s(x)\cdot c$ and the extension $\esatta 1 \Gamma \i G \p H 1$ is isomorphic to $\esatta 1 \Gamma {} {H\times_h \Gamma} {} H 1$, where $H\times_h \Gamma \to H$ is the projection $pr_{1}$ on the first coordinate. 
Note that if the 2-cocycle $h$ is definable in $(\Gamma,\M)$, then $H\times_h\Gamma$ is definable in $(\Gamma,\M)$, and this gives an interpretation of $\str G \pi \M$ in $(\Gamma,\M)$. We have thus proved:

\bprop \label{coc} If $\esatta 1 \Gamma \i G \p H 1$ admits a section $s\colon H \to G$ such that the corresponding cocycle $h\colon H\times H \to \Gamma$ is definable in  $(\Gamma,\M)$, then $\str G \pi \M$ is naturally interpretable in $(\Gamma,\M)$. Any parameters required for the interpretation are those needed to define $h$ in $(\Gamma,\M)$.
\eprop

An important situation where the hypothesis (and so also the conclusion) of Proposition 3.1 holds is when $M$ is an $o$-minimal expansion of $\rf$, $H$ is a connected Lie group definable in $M$ and $\pi:G \to H$ is the universal cover of $H$. This will be discussed in detail at the beginning of subsection 3.2.
A similar situation holds when $M$ is an $o$-minimal expansion of an arbitrary real closed field, $H$ a definably connected group definable in $M$, and $G$ is the $o$-minimal universal cover of $H$, which is by definition a locally definable group in $M$. 

\vspace{2mm}
\noindent
Finally in this subsection we point out that interpretability as in Proposition 3.1 for one model yields a stable embeddability result at the level of theories.

\bprop \label{uniform} Suppose that $\esatta 1 \Gamma \i G \p H 1$ admits a section $s\colon H \to G$ such that the corresponding cocycle $h\colon H\times H \to \Gamma$ is definable in  $(\Gamma,\M)$. By adding constants for suitable elements  from $\Gamma \leq G$ and $M$, assume that $h$ is $\emptyset$-definable in $(\Gamma,\M)$. 
\newline
Then for any $\str {G'} {\p'} {\M'} \equiv \str G \p \M$, $M'$ is stably embedded in $\str {G'} {\p'} {\M'}$.
\eprop 
\bp It is enough to prove the conclusion when $\str {G'} {\p'} {\M'}$ is saturated. There is no harm in assuming that our languages and theories are countable. We fix some big cardinal $\kappa$ such that any countable theory has a 
(necessarily unique) saturated model of cardinality $\kappa$. As mentioned in the introduction, the fact that  
$(H\times_{h}\Gamma,pr_{1}, M)$ is a model of $Th(G,\pi,M)$ is contained in $Th(\Gamma,\M)$. Let $(\Gamma',\M')$ be a saturated model of cardinality $\kappa$ of $Th(\Gamma,\M)$. So the corresponding  $(H'\times_{h'}\Gamma', pr_{1},M')$ is also saturated, of cardinality $\kappa$ and elementarily equivalent to $\str G \p \M$. So it suffices to show that 
$M'$ is stably embedded in $(H'\times_{h'}\Gamma', pr_{1},M')$. And for that it suffices to prove that $M'$ is stably 
embedded in $(\Gamma',\M')$. This is clear, but we give a proof anyway. First if $a,m$ are tuples from $\Gamma',M'$ 
respectively, then $tp(a,m)$ in $(\Gamma',\M')$ is determined by $tp(a)$ in $\Gamma'$ and $tp(m)$ in $M'$. By 
compactness any formula $\phi(x,y)$ in the language of $(\Gamma',M')$, (where $x,y$ are tuples of variables ranging 
over tuples of the appropriate length from $\Gamma'$, $\M'$, respectively) is equivalent, modulo $Th((\Gamma',\M')$ 
to a finite disjunction of formulas of the form 
$\theta(x)\wedge\psi(y)$ where $\theta$ is in the language of $\Gamma'$ and $\psi$ in the language of $M'$. Now let $\phi(x,y)$ be in the language of $(\Gamma',\M')$, let us fix some tuple $a$ from $\Gamma'$, and we have to show that the set $Y$ of tuples $m$ from $M'$ such that $(\Gamma',M')\models \phi(a,m)$, is definable (possibly with parameters) in the structure $M'$. Suppose $\phi(x,y)$ is equivalent, as above, to the finite disjunction 
$\bigvee_{i\in I}\theta_{i}(x)\wedge \psi_{i}(y)$. Let $I_{0}\subseteq I$ be the set of $i$ such that 
$\Gamma'\models \theta_{i}(a)$. Then $Y$ is definable in $M'$, without parameters, by the formula 
$\bigvee_{i\in I_{0}}\psi_{i}(y)$. 
\ep

Warning: It may happen that some models $\str {G'}{\pi'}{\M'}$ of $Th(\str G \pi \M)$ (even with $\M'=\M$ and $\Gamma'=\Gamma$) are not interpretable in a model of $Th(\Gamma,\M)$.
The problem is that there could be no section $s\colon H\to G'$ definable in $\str {G'}{\pi'}{\M'}$ inducing the given 2-cocycle $h$. An example is the exotic extension of $\R/\Z$ by $\Z$ given in Theorem \ref{modZ}.

\subsection{Relative $L_{\omega_1,\omega}$-categoricity}
When we talk about categoricity in the language $L_{\omega_1,\omega}$, we are by convention working over a countable language $L$, and we are usually interested in forms of categoricity (in a given cardinal, relative to the isomorphism type of part of the structure,  absolute,...)
of a {\em single} $L_{\omega_1,\omega}$ sentence $\sigma$. So for example any countable $L$-structure $M$ is  
$L_{\omega_1,\omega}$-$\omega$-categorical in the sense that there is a single $L_{\omega_1,\omega}$-sentence $\sigma$ whose unique countable model is $M$. If in addition all elements of $M$ are named by constants then $M$ is absolutely $L_{\omega_1,\omega}$-categorical: there is a $L_{\omega_1,\omega}$-sentence whose unique model is $M$. 

Our main result is Theorem 3.4 below which concerns categoricity in $L_{\omega_1,\omega}$, relative to the isomorphism type of the field, in a strong sense.

The ``standard" model we will be concerned with is  $(G,\pi,M)$ where for convenience $M$ is just $\rf$, and $\pi:G\to H$ is the universal cover of some connected Lie group $H$ which is definable in $M$ (and we add constants for parameters over which $H$ is defined). So $G$ is the universal cover of $H$ as a topological (equivalently Lie) group, but in the structure $(G,\pi, M)$ on the face of it we only see the group structure on $G$. It was proved in \cite{HrPePi:08b} that Proposition 3.1 applies to this situation.
This is the content of section 8.1 of \cite{HrPePi:08b} and specifically of Theorem 8.5 there and its proof which we give a brief summary of. The key point (depending on results of Edmundo and Eleftheriou) is that the universal cover $\pi:G\to H$ of $H$ can be realized as a topological group, $\w H$ say, which is ``locally definable" (or $\bigvee$-definable) in $M$.  In other words there is a locally definable group $\w H$ in $M$  (equipped via $o$-minimality with topological so Lie group structure) and a locally definable surjective homomorphism ${\w \pi}:{\w H} \to H$, and there is an isomorphism of topological groups $f:G\to \w H$  which makes everything commute. Let $\Gamma_{1} < \w H$ be $f(\Gamma)$, so a discrete group also locally definable in $M$.
Now $\w H$ being locally definable in $M$, there is (using the existence of Skolem functions) a section $s:H\to \w H$ definable in $M$. 
The corresponding cocycle $h:H\times H \to \Gamma_{1}$ has finite image $\Gamma_{0} < \Gamma_{1}$. So the map $h:H\times H \to \Gamma_{0}$ is a partition of $H\times H$  which is again definable in $M$. Hence the locally definable group $\w H$ is isomorphic to (so can be identified with) $H\times_{h}\Gamma_{1}$ which is definable in 
$(\Gamma_{1},\M)$. Identifying  $\Gamma$ with $\Gamma_{1}$ this gives an isomorphism $f: G \to H\times_{h}\Gamma$, and gives the interpretation of $(G,\pi,M)$ in $(\Gamma, \M)$. Note that $\Gamma$ is $\emptyset$-definable in $(G,\pi,M)$, but there is no reason to believe that $f$ is definable in $(G,\pi,M)$ (in fact in general it is not). Note that $f$ commutes with the canonical projection maps  $\pi:G\to H$ and $pr_{1}: H\times_{h}\Gamma  \to H$. In this situation, the key new lemma needed for the relative $L_{\omega_1,\omega}$-categoricity statement is:

\bl \label{infsection} The isomorphism $f: G \to H\times_{h}\Gamma$ is $L_{\omega_1,\omega}$-definable without parameters in the structure
$(G,\pi,M)$. Namely there is an $L_{\omega_1,\omega}$-formula $\chi(x,y)$ such that for any $a\in G$ and $b\in H\times_{h}\Gamma$, $f(a) = b$ iff  $(G,\pi,M)\models \chi(a,b)$. 
\el 

Let us note in passing that using the isomorphism $f: G \to H\times_{h}\Gamma$ we can define a section $s\colon H\to G$, and conversely given a section $s$ we can define $f$ (with $h$ the cocycle induced by the section). This observation will be used later. Let us now prove the Lemma.

\bp 
Let $\Lambda$ range over finite index subgroups of  $\Gamma$ which are $\emptyset$-definable in $(\Gamma,+)$. As 
$\Gamma$ is a finitely generated abelian group, the intersection of all such $\Lambda$ is $0$. For each $\Lambda$, $f$ induces an isomorphism $f_{\Lambda}:G/\Lambda \to  (H\times_{h}\Gamma)/\Lambda$. 
\begin{claim}  Each $f_{\Lambda}$  is $\emptyset$-definable in $(G,\pi,M)$. 
\end{claim}
\bp  Let $n$ be the index of $\Lambda$ in $\Gamma$. We have surjective homomorphisms from $G/\Lambda$ to $H$ and  $(H\times_{h}\Gamma)/\Lambda$ to $H$, induced by $\pi, pr_{1}$ which we will still call $\pi, pr_{1}$. Moreover 
$f_{\Lambda}$ commutes with these projections.
Note that $\pi: G/\Lambda \to H$ is, among other things, a finite topological cover (of connected Lie groups). By Proposition 2.2 for example, the group $G/\Lambda$ is  interpretable in $M$, so in particular it can be considered as a definably connected group definable in $M$. By 
\cite[Lemma 8.1]{HrPePi:08b}, for some $k$, every element of $G/\Lambda$ can be written as a product of $k$ $nth$ powers.   Let $y\in G/\Lambda$, and write $y = y_{1}^{n}\cdot \ldots \cdot y_{k}^{n}$. For $i=1,\ldots, k$, let $z_{i}\in (H\times_{h}\Gamma)/\Lambda$  such that $pr_{1}(z_{i}) = \pi(y_{i})$. As $n = |\ker(\pi)| = |\ker(pr_{1})|$ and the extensions are central, $z_{i}^{n}$ depends only on $\pi(y_{i})$. Clearly $f_{\Lambda}(y) = z_{1}^{n}\cdot \ldots \cdot z_{k}^{n}$. 
\ep

We now complete the proof of Lemma 3.3. Let $\chi(x,y)$ be the conjunction of  $f_{\Lambda}(x/\Lambda) = y/\Lambda$, where $\Lambda$ ranges over the (countable) family of $\emptyset$-definable subgroups of $\Gamma$ of finite index. 
\ep 
Let us now prove the main result of this section. Recall that  
$\pi:G\to H$ is the universal cover of some connected Lie group $H$ which is definable in $M=\rf$ 

\bt \label{categorical} Consider the structure $(G,\pi,M)$ equipped also with constants for generators of $\Gamma$, and let $L$ be its language. Then 
\newline
(i) there is a single  $L_{\omega_1,\omega}$-sentence $\sigma$ which is true of $(G,\pi, M)$ and such that whenever  $(G_{1},\pi_{1},M_{1})$, $(G_{2},\pi_{2},M_{2})$ are models of $\sigma$, then ANY isomorphism  between $M_{1}$ and $M_{2}$ extends (uniquely) to an isomorphism between $(G_{1},\pi_{1},M_{1})$ and $(G_{2},\pi_{2},M_{2})$.
\newline
(ii) Let $T'$ be the set of {\em all}  $L_{\omega_1,\omega}$ sentences true in $(G,\pi, M)$. Then $(G,\pi,M)$ is the unique model of $T'$. 
\et 
\bp  
(i)  Note first that if  $(G',\pi',M')$ is (first order) elementarily equivalent to $(G,\pi,M)$, then we have the group $H'\times_{h'}\Gamma'$ $\emptyset$-definable in $(G',\pi',\M')$.  Moreover if $(G'',\pi'',\M'')$ is elementarily equivalent to $(G,\pi,M)$ and $c$ is an isomorphism between $\Gamma'$ and $\Gamma''$ and $d$ an isomorphism between $M'$ and $M''$  (fixing the appropriate constants) then $(c,d)$ determines an isomorphism  between the structures $(\Gamma',\M')$ and $(\Gamma'',\M'')$ and in particular an isomorphism which we also call $(c,d)$ between the groups  $H'\times_{h'}\Gamma'$ and
$H''\times_{h''}\Gamma''$.  So we let $\sigma$ be the conjunction of  (i) the first order theory of $(G,\pi,M)$, (ii) a sentence fixing the isomorphism type of $\Gamma$ together with its generators, (iii) a sentence expressing that 
$\chi(x,y)$ is an isomorphism between $G$ and $H\times_{h}\Gamma$. Hence if $(G',\pi',M')$ and $(G'',\pi'',\M'')$ are 
models of $\sigma$, then we already have an isomorphism $c:\Gamma' \to \Gamma''$, and if $d: \M' \to \M''$ is an 
isomorphism, then $\chi((G'',\pi'',\M''))^{-1}\circ (c,d) \circ \chi((G',\pi',\M'))$ is an isomorphism between $G'$ 
and $G''$ which together with $d$  gives an isomorphism between $(G',\pi',\M')$ and $(G'',\pi'',\M'')$, as required.  
(In the above the notation $\chi((G',\pi',\M'))$ is supposed to denote the interpretation of $\chi$ in $(G',\pi',\M')$.)

\vspace{2mm}
\noindent
(ii) Note that $M = \rf$  (with finitely many additional constants) can be characterised up to isomorphism by the collection of {\em all} $L_{\omega_1,\omega}$-sentences true in it. (In addition to the first order theory, say that the rationals are dense and that every Dedekind cut in the rationals is realized.)  Together with $\sigma$ from part (i), this collection of infinitary sentences has a unique model. 
\ep 

We end with a couple of remarks. First the formal statement of Theorem 3.4(ii) is not in itself surprising: given a 
structure $N$ of cardinality the continuum, it may happen (as with the reals) that every element of $N$ is     
  ``$L_{\omega_1,\omega}$-definable", in which case $N$ is clearly the unique model of its full $L_{\omega_1,\omega}$ theory.
  So it is Theorem 3.4(i) which we consider the main point.  We can then ask if the analogous statement holds in the complex case. We consider now the Zilber structure  $((\C,+), exp, (\C,+,\cdot))$  (where there is no connection between the two sorts other than $exp$), with say an additional constants for the generator of $\ker(exp)$. The question is whether there is a single $L_{\omega_1,\omega}$-sentence $\sigma$ of the appropriate language, such that given models
  $(V,\pi,K)$ and $(V',\pi',K')$ of $\sigma$ any isomorphism between the field sorts $K$ and $K'$ lifts to an 
isomorphism between the full structures. We conjecture that the answer is in fact NO. The categoricity statement that 
Zilber actually proves is that if $(V,\pi,K)$ and $(V',\pi',K')$ are models of the first order theory of the standard 
model, and both have standard kernels  (and this is the only property that needs to be expressed by an infinitary sentence) then if $K, K'$ are isomorphic (equivalently 
have the same cardinality) then  the two structures $(V,\pi,K)$ and $(V',\pi',K')$ are isomorphic by an isomorphism $f$
 say. But we cannot demand that $f$ lifts any given isomorphism between $K$ and $K'$. 
  
 In the $o$-minimal setting, we have considerably more rigidity in the field sort. For example there is only one isomorphism between two copies of the real field, and relative categoricity in the sense of Zilber, will fail. The rest of the paper is devoted towards clarifying and expanding on some of these issues.

\subsection{Inverse limits}
\label{dense}

\bd Let $H$ be a real Lie group, say connected, let $$\ex \Gamma {\w H} H$$ be the universal cover of $H$, and let $\widehat H$ be the inverse limit of all the finite covering homomorphisms of $H$. Up to isomorphism each finite cover of $H$ is a quotient of the universal cover. So $$\widehat H = \varprojlim_{\Lambda \in \ca J} {\w H}/\Lambda$$ where $\ca J$ is the family of all subgroups of $\Gamma$ of finite index. We can and will identify $\widehat H$ with the subgroup of $\Pi_{\Lambda \in \ca J} {\w H}/\Lambda$ consisting of those elements $(x_\Lambda)_{\Lambda\in \ca J}$ of the product such that whenever $\Lambda \subset \Lambda'$ the element $x_\Lambda$ is mapped to $x_{\Lambda'}$ under the natural homomorphism ${\w H}/\Lambda \to {\w H}/\Lambda'$. \ed

With the notation above we say that a subgroup $G$ of $\widehat H$ is {\em dense} if for any 
$x = (x_\Lambda)_{\Lambda\in \ca J}$ in ${\widehat H}$ and finite subset $F$ of $\ca J$, there is $y\in G$ such that 
$y_{\Lambda} = x_{\Lambda}$ for  $\Lambda \in F$. This is clearly equivalent to the natural maps (in fact homomorphisms) from $G$ to 
${\w H}/\Lambda$ being surjective for all  $\Lambda \in \ca J$. So we obtain immediately:

\bprop \label{density} Let $H$ be a connected real Lie group and let $$\ex \Gamma {\w H} H$$ be the universal cover of $H$. The natural embedding $f\colon \w H \to \widehat H$, induced by the morphisms from $\w H$ to the finite covers of $H$, has dense image in $\widehat H$. \eprop

\bexa 
Each finite cover of $\R/\Z$ is of the form $\R/\Z\stackrel k \longrightarrow \R/\Z$ where ``$k$'' is multiplication by the positive integer $k$. Let $$\widehat {\R/\Z} \subset \Pi_{k\in \NN^*} \R/\Z$$ be the inverse limit of these finite covers. The elements of $\widehat {\R/\Z}$ are the sequences $(h_i : i>0) \in \Pi_{i\in \NN^*}{\R/\Z}$ satisfying the compatibility conditions $kh_{ik} = h_i \in \R/\Z$ for all positive integers $i,k$. 
\eexa

\bprop \label{good-sequences} Let $$\esatta 0 \Z \i {(\R,+)} \modZ {\R/\Z} 0$$ be the universal cover of the circle group.
There is a natural embedding $f\colon \R\to \widehat {\R/\Z}  \subset \Pi_{n\in \NN^*}{\R/\Z}$ sending $x$ to $(\bmod_{\Z} (x/n) : n>0)$. The image of $\R$ is dense in $\widehat {\R/\Z}$ and consists exactly of those sequences $(h_i : i>0)\in \widehat {\R/\Z}$ such that $h_i\to 0$ in $\R/\Z$ for $i\to \infty$.
\eprop

\bp The density is stated already in Proposition \ref{density}. We prove the second part. One direction is obvious: if $x\in \R$ then $x/n \to 0$ in $\R$, so the images in $\R/\Z$ also tend to $0$.  Conversely, let $(h_i : i>0)\in \widehat {\R/\Z}$ and suppose that $h_i\to 0$. We must find $x\in \R$ such that $\pi(x/n) = h_n$ for every $n$. 
Since $h_i\to 0$ for $i\to \infty$, in particular $h_{2^n}\to 0$ for $n\to \infty$. So there is some $n_0\in \NN$ such that for all $n\geq n_0$ we have $h_{2^n} \in V$ where $V$ is the image of $(-1/4,1/4)$ under the natural projection $\R \to \R/\Z$. Note that for each $a\in V$ the equation $2x=a$ has a unique solution $x$ in $V$, which we call $\frac 1 2 a$ (there is of course a second solution in $\R/\Z$ lying outside of $V$). So for all $n \geq n_0$ we have $h_{2^{n+1}} = \frac 1 2 h_{2^n}$. It then clearly follows that there is a real number $x\in \R$ such that $\modZ(x/m) = h_m$ for all powers of two $m=2^n$ bigger than $2^{n_0}$. Now consider the sequence $a_i:=\modZ(x/i)$, $i\in \NN^*$. It suffices to show that $a_i = h_i$ for every $i \in \NN^*$. So consider the difference $b_i = a_i-h_i$. Then $b_i$ tends to zero for $i\to \infty$, and $b_i = 0$ for all powers $i$ of $2$ bigger than $2^{n_0}$. So it suffices to prove the following claim. 

Claim: Let $(b_i:i>0)\in \widehat{\R/\Z}$ satisfy the following two conditions: (i) $b_i\to 0$ for $i\to \infty$; (ii) $b_i = 0$ for infinitely many values of $i$. Then $b_i = 0$ for all $i$. 

To prove the claim let $(n_k)_{k\in \NN^*}$ be an infinite sequence with $b_{n_k} = 0$ for all $k\in \NN^*$. Given $m$ and $k$, consider $b_{m\cdot n_k}$. Since
$b_{n_k}=0$,  $b_{m\cdot n_k}$ should be an $m$-th torsion element. As $k$ tends to
infinity, $b_{m\cdot n_k}$ should tend to zero, so the only $m$-th root it can be is the
zero element. Hence $b_{m\cdot n_k}$ is eventually zero for large $k$. But then
$b_m=n_k\cdot b_{m \cdot n_k}=0$. 
\ep

\subsection{The universal cover of $\R/\Z$}
We study in more detail the model  theory of the universal cover $\esatta 0 \Z {} \R \modZ {\R/\Z} 0$ of the circle group. We do this both (i) in the ``pure group" language where both $\R$ and $\R/\Z$ are equipped only with their (additive) group structures (and of course we have a symbol for the covering map $\bmod_{\Z}$) and (ii) 
in the richer language in which the group $\R/\Z$ is viewed as $([0,1), +(mod\Z))$ equipped with predicates for all sets of $n$-tuples $\emptyset$-definable in $\rf$. In fact in case (ii) there is no harm in working, as earlier, with the structure $((\R,+),\modZ,\rf)$  where of course only  $\bmod_{\Z}$  links the two sorts. (We can actually work in any sub-language in which the topology of $\R/\Z$ is definable, so for instance we can replace $\rf$ by $(\R, +, <)$ and work with the structure $((\R,+),\modZ,(\R,+,<))$.)
 
If we work in context (i), it turns out that the only extension of $\R/\Z$ by $\Z$ which is a model of the first order theory of $\esatta 0 \Z {} \R \modZ {\R/\Z} 0$ is (up to isomorphism) the standard model $\str {(\R,+)}\modZ {\R/\Z}$ (this also follows from Zilber's results). However in context (ii) 
there are at least $2^{\aleph_0}$ pairwise non-isomorphic models. This implies in particular that in the structure $((\R,+),\modZ,\rf)$ we cannot define a section $s \colon \R/\Z =[0,1) \to \R$ in a first order way, for otherwise we could use $s$ to obtain an isomorphism $f$ as in Lemma \ref{infsection} and then obtain a first-order categoricity result reasoning as in Theorem \ref{categorical}. So roughly speaking our results say that we cannot lift the topology of $\R/\Z$ to $\R$ in a first order way, although (using the $L_{\omega_1,\omega}$-section) we can do it by an infinitary formula. Let us begin with context (i). We need:

\bl \label{hull} Let $V$ be an abelian divisible torsion free group (hence a $\Q$-vector space). Let $\Gamma$ be an infinite subgroup of $V$ with $|\Gamma|<|V|$. Then the isomorphism type of $(V,+,\Gamma)$ is determined by the isomorphism type of $(\Gamma,+)$ and the isomorphism type of $(V/\Gamma, +)$. \el 
\bp Suppose that $(V_1,+,\Gamma_1)$  and $(V_2,+,\Gamma_2)$ have $\Gamma_1\cong \Gamma_2$ and $V_1/\Gamma_1\cong V_2/\Gamma_2$. For $i=1,2$, the $\Q$-vector subspace $\langle \Gamma_i\rangle$ of $V_i$ generated by $\Gamma_i$ is isomorphic to $\Gamma_i \otimes_\Z \Q$ and therefore its isomorphism type is determined by the isomophism type of $\Gamma_i$. So any isomorphism $f\colon \Gamma_1\to \Gamma_2$ extends to an isomorphism $\ov f \colon \langle \Gamma_1\rangle \to \langle \Gamma_2 \rangle$. By our assumptions $|V_1/\Gamma_1| = |V_2/\Gamma_2|$. Since $\aleph_0 \leq |\Gamma_i|<|V_i|$, this is equivalent to $\dim(V_1/\langle \Gamma_1 \rangle) = \dim(V_2/\langle \Gamma_2 \rangle)$. So $\ov f$ extends to an isomorphism $(V_1, + ,\Gamma_1) \cong (V_2, + , \Gamma_2)$. \ep

In the above proof note that the isomorphism $(V_1, + ,\Gamma_1) \cong (V_2, + , \Gamma_2)$ can be chosen to extend any given isomorphism $f\colon \Gamma_1\to \Gamma_2$.

\bc \label{pure} Let $(G,+)$ be an abelian divisible torsion free group and consider an exact sequence 
$E$ of the form $$\esatta 0 {\Z} \iota  {(G,+)} \pi {\R/\Z} 0.$$  In the pure group language $E$ is isomorphic to the universal cover of the circle group
$$\esatta 0 {\Z} {}  {(\R,+)} \modZ {\R/\Z} 0.$$ Indeed there is such an isomorphism which fixes pointwise the kernel $\Z$ and permutes $\R/\Z$ by a group automorphism. \ec

\bc \label{compatible} With $G$ as above, let $\hat \pi \colon G \to \widehat {\R/\Z}$ be the map sending $g$ to $(\pi(g/n))_{n\in \NN^*}$. Then $\hat \pi$ has dense image. 
\ec 
\bp By Proposition \ref{good-sequences} and Corollary \ref{pure}. \ep

Let us now consider context (ii). To repeat, our ``standard model'' is 
\newline
$((\R,+), \modZ, \rf)$. In addition to the field structure on $\rf$, the additive group structure on the first copy $\R$ of the reals, and  the covering map 
$\bmod_{\Z}$, it will be convenient to adjoin a constant symbol for a given generator $1_\Z$ of the kernel $\Z$ of 
$\bmod_{\Z}$. We will call this language $L$. 

\bt \label{modZ} Let $(G,+)$ be an abelian divisible torsion free group and consider an exact sequence 
$E$ of the form $$\esatta 0 {\Z} \iota  {(G,+)} \pi {\R/\Z} 0$$ viewed naturally as an $L$-structure $(G,\pi,\rf)$.  Then:  
\begin{enumerate}

\item $\str G \p \rf \equiv \str {(\R,+)} \modZ \rf$ in the language $L$ if and only  if for every positive $n\in \NN$, $\pi(\iota(1_\Z)/n) = \modZ(1/n)$ (which equals $1/n$ in the field sort under the identification of $\R/\Z$ with $[0,1)$).  
\item The isomorphism type of $(G,\pi,\rf)$ is given precisely by the image of the map $\hat \pi\colon G \to \widehat {\R/\Z}$, $x\mapsto (\pi(x/n))_{n\in \NN^*}$, and there are at least $2^{\aleph_0}$ possibilities.
\item The structure $(G,\pi,\rf)$ is isomorphic (in the language $L$) to the standard model  $((\R,+), \modZ, \rf)$ if and only if $(G,+)$ can be equipped with the structure of a topological group making $\pi$ a covering homomorphism. 
\end{enumerate}
\et

\bp (1) Note first that the fact that $\modZ (1_{\Z}/n) = (1/n)$ for all $n > 0$  is part of $Th((\R,+), \bmod_{\Z}, \rf)$ (as  $1_{\Z}$ is named by a constant, and  $1/n$ in the field sort is also named by a constant). Hence we have left to right.

To prove the converse we can take two saturated elementary extensions of the respective structures, and show they are isomorphic. In fact it is no more work to give axioms for what will be the common theory. 
So we will define a theory $T$, with sorts $\Gamma,G,M$, such that both $\str G \pi \rf$ and $\str {(\R,+)} \modZ \rf$ are models of $T$ (with a constant for $1_\Z$), and then we will prove the completeness of $T$  by showing that two saturated models (of the same cardinality) are isomorphic. 
The theory $T$ says the following: 
\begin{enumerate}
\item $M\equiv \rf$, 
\item $H := \pi(G) = [0,1)^M$ with addition modulo $1$, 
\item $\esatta 0 \Gamma \i G \p H  0$ is an exact sequence of abelian groups, 
\item $G$ is divisible and torsion free, 
\item $1_\Gamma$ is an element of $\Gamma$ such that $h_n:= \pi(\i(1_\Gamma/n))$ is $1/n \in H\subset M$. 
\end{enumerate}
To prove the completeness of $T$ fix a saturated model $N = \structure \Gamma \i G \p \M$ of cardinality $\kappa > \aleph_0$, let $H = \pi(G)$ and let $\widehat H$ be the inverse limit of the system of maps $k \colon H \to H, x \mapsto kx$. The elements of $\widehat H$ are the sequences $(a_n \mid n >0)$ with $a_n\in H$ such that $ka_{nk} = a_n$ for all positive integers $k,n$. Let $\widehat \pi \colon G \to \widehat H$ be the map $g \mapsto (\pi(g/n):n >0)$. 
By Corollary 3.10, the reduct $N_{0}$ of $N$ to the ``pure group language" is isomorphic to a saturated elementary extension of the standard model $\str {(\R,+)}\modZ {\R/\Z}$ discussed earlier. So by Corollary 3.11 and compactness (i.e. saturation), 
the map $\widehat \pi \colon G \to \widehat H$  is surjective. 
Let $$V:= \ker(\widehat \pi).$$ Then $V$ is a $\Q$-vector space of dimension $\kappa$ included in $\ker(\pi)$. (Note that in the standard model $V=0$.) Since $G$ is a $\Q$-vector space we can write $$G = V \oplus W$$ for some $\Q$-vector space $W< G$. The subgroup $W$ is not unique, but for any choice of $W$ we have $$W \cong \widehat H$$ via $\widehat \pi$. Similarly, given another saturated model $N_{1} = \structure {\Gamma_1} {\iota_1} {G_1}  {\pi_1} {\M_1}$ of cardinality $\kappa$, we obtain $G_1 = V_1 \oplus W_1$ as above. Since the theory of $T$ includes the complete theory of $\rf$, we have $\M \cong \M_1$ (hence $H\cong H_1$). To prove the proposition we can as well assume that $\M=\M_1$. To simplify notations we can also assume that $\iota, \iota_1$ are the inclusion maps. Since $V,V_1$ are $\Q$-vector spaces of the same dimension, we have $V\cong V_1$. Moreover $W \cong \widehat H \cong W_1$ by composing $\widehat \pi$ with the inverse of $\widehat {\pi_1}$. We conclude that $G \cong G_1$ as coverings of $H$ (namely the isomorphism obtained in this way commutes with $\pi,\pi_1$). Note that the proof so far yields completeness of the ``reduct" $T_{0}$ of $T$ to the language which omits the constant symbol for a generator of $\Z$.  

But, working in the language $L$, we must still show that, for a suitable choice of $W,W_1$, the isomorphism $G\cong G_1$ constructed above sends $1_\Gamma$ to $1_{\Gamma_1}$ (clearly it sends $\ker(\pi)$ to $\ker(\pi_1)$, but this is not enough). To this aim let $\Z1_\Gamma$ be the subgroup of $G$ generated by $1_\Gamma$. Note that $\Z1_\Gamma \subset \ker(\pi)$ (the equality holds in the standard model). Now let $\Q1_\Gamma$ be the divisible hull of $\Z1_\Gamma$ and observe that $V \cap \Q1_\Gamma = 0$ (if not $\pi(1_\Gamma/n) = 0$ for infinitely many $n$, contradicting the axiom $\pi(1_G/n) = h_n$). It then follows that we can write $G = V \oplus W$ for some $\Q$-vector space $W<G$ containing $\Q1_\Gamma$. In particular $1_\Gamma \in W$. In $N_{1}$ we can similarly write $G_1 = V_1 \oplus W_1$ with $1_{\Gamma_1}\in W_1$. As above we have an isomorphism $G \cong G_1$ commuting with $\pi,\pi_1$ and sending $V$ to $V_1$ and $W$ to $W_1$. This isomorphism must send $1_\Gamma$ to some element $x\in W_1$ with $\pi_1(x/n) = h_n$ for every $n$. But $1_{\Gamma_1}$ is the unique such $x$. So $1_\Gamma$ goes to $1_{\Gamma_1}$ and we are done. Namely we have exhibited an isomorphism between the $L$-structures $N$ and $N_{1}$. 

\medskip 
\noindent
(2) Suppose that $E'$ is another exact sequence  which is isomorphic to $E$ as an $L$-structure, via isomorphism $f$ say. Note that $f$ must be the identity on $\rf$  (also on $\Z$ assuming that the same generator is named by the constant in the two structures). 
Hence for each  $a\in G$, and $n > 0$, $f(\pi(a/n)) =  \pi'(f(a)/n)$. This says precisely that $\hat \pi (G)$ and ${\widehat \pi'}(G')$ have the same image in ${\widehat \R/\Z}$. 

Conversely if $E'$ is another exact sequence  with kernel $\Z$, and image $\R/\Z$ and $G'$ divisible torsion-free then the map ${\widehat \pi'}:  G' \to \widehat{\R/\Z}$ is well-defined, and is an embedding. 
So  $(G',\pi',\rf)$ is isomorphic to $({\widehat \pi'}(G'), pr_1, \rf)$  where  $pr_1$ is projection on the first coordinate. 

We have shown the first part of (2). For the rest of it, we construct continuum many suitable subgroups of $\widehat{\R/\Z}$. Our standard model is  $((\R,+), \modZ, \rf)$ with an additional constant for $1$ in the first sort $\R$.  We have the embedding 
$\widehat{\modZ}: \R \to \widehat{\R/\Z}$. In particular we have $\widehat {\modZ}(\Z) \leq \widehat{\R/\Z}$.
Let $H$ denote the subgroup of $\widehat{\R/\Z}$ consisting of sequences ${\bar h} = (h_{n})_{n}$ with $h_{1} = 0$. Then $H$ is of cardinality continuum, whereas $\widehat{\modZ}(\Z)\leq H$ is countable. 
Let $\{{\bar h}^{\alpha}:\alpha < 2^{\aleph_{0}}\}$ be a set of continuuum many elements of $H$ which are in different cosets modulo $\widehat{\modZ}(\Z)$. 
Let $g\in \R \setminus \Q$ and let ${\bar v} = \hat\pi(g) \in \widehat{\R/\Z}$. For each $\alpha$, let ${\bar w}^{\alpha} =  {\bar v} + {\bar h}^{\alpha}  \in \widehat{\R/\Z}$. 
\newline
{\em Claim.} For each $\alpha$ there is $G^{\alpha}\leq \widehat{\R/\Z}$ such that ${\bar w}^{\alpha}\in G^{\alpha}$, 
$G^{\alpha} \cap H = \widehat{\modZ}(\Z)$ and  $(G^{\alpha}, pr_{1}, \rf)$ (with the interpretation of the distinguished constant $1$ as  $\modZ(1_{\Z})$) is elementarily equivalent to the standard model (equivalently is a model of the theory $T$ described above). 
\bp Note that NO element of the sequence ${\bar w}^{\alpha}$ is a torsion element of $\R/\Z$. Hence the set of integer multiples of elements of ${\bar w}^{\alpha}$ is a $1$-dimensional $\Q$-vector space, $W^{\alpha}$ say,  which is disjoint over $0$ from the torsion subgroup $T$ of $\R/\Z$. Hence $W^{\alpha}$ extends to a $\Q$-vector space
$W'^{\alpha}$ such that  $T \oplus W'^{\alpha} = \R/\Z$. Let $G^{\alpha} = 
\widehat \modZ (\Q) \oplus \{(v/n)_{n\in \NN^*}: v\in W'^{\alpha} \} \cong \Q \oplus W'^{\alpha}$. 
\ep

Note that for $\alpha \neq \beta$, ${\bar w}^{\beta}\notin G^{\alpha}$. For otherwise 
${\bar w}^{\beta} - {\bar w}^{\alpha} = {\bar h}^{\beta} - {\bar h}^{\alpha}  \in G^{\alpha}$. But 
${\bar h}^{\beta} - {\bar h}^{\alpha} \in H\setminus \widehat \modZ (\Z)$, contradicting the properties of $G^{\alpha}$ in the Claim above.  Hence  $G^{\alpha} \neq G^{\beta}$ for $\alpha \neq \beta$. This yields (2).

\vspace{2mm}
\noindent
(3)  This is really a straightforward and easy group-theoretic/topological remark, which is well-known:  Namely, IF that $(G,+)$ is an abelian divisible topological group, equipped with a covering homomorphism $\pi$ onto the topological group $\R/\Z$ and with kernel $\Gamma$ abstractly isomorphic to $\Z$,  THEN $\pi: G\to \R/\Z$  IS the universal cover. 
\newline
Let us note in passing  that $G$ has the structure of a $1$-dimensional Lie group. Also $\Gamma$ is a discrete subgroup of $G$, and  $\pi:G \to \R/\Z$ is the quotient map (as a map of topological, or Lie groups). 
\newline
By the universal property of the universal cover $\modZ \colon \R \to \R/\Z$ there is a (unique) continuous homomorphism 
$f:\R \to G$ such that $\pi\circ f = \modZ$. Note that $\ker(f) \leq \Z$. So $f(\R)$ is a connected $1$-dimensional Lie 
group, which is torsion-free, hence homeomorphic to $\R$, in particular a $\Q$-vector space. But then $G/f(\R)$ is 
both a $\Q$-vector space as well as a quotient of $\Gamma$, so has to be trivial. Namely $f:\R \to G$ is surjective. Now a
$\ker(f)$ is a subgroup of  $\ker(\modZ) = \Z$, so is either trivial, or a finite index subgroup of $\Z$. In the latter 
case $f(\R)$ is compact, impossible, so $\ker(f)$ is trivial, and $f$ is a homeomorphism. 

\ep

\vspace{2mm}
\noindent
Let us finish this subsection by remarking that the non categoricity (relative to the kernel $\Z$, and to the real field $\rf$) statements in Theorem 3.12 also hold with the reals replaced by an arbitrary (even saturated) real closed field, even though such a field may have automorphisms.

\subsection{Questions and final comments}
We end with a couple of problems. 
\begin{problem} Describe the groups definable in (a) the $2$-sorted structure  $((\Z,+),\rf)$, or more generally in 
(b) $((\Z,+), M)$ for any $o$-minimal expansion $M$ of $\rf$, or more generally in (c) a structure $((A,+),M)$ where
$(A,+) \equiv (\Z,+)$ and $M$ is $o$-minimal expansion of an arbitrary real closed field.
\end{problem}

\noindent
Comment: We have seen that universal covers of simple Lie groups are examples of groups definable in the structure (a).

\begin{problem}
Let $\M$ be an $o$-minimal expansion of  $\rf$, and
Let $H$ be a definable, definably connected group in $\M$. 
Let $$\exact \Gamma \i G \p H$$ be an exact sequence of groups, with $\i (\Gamma) < Z(G)$ and $\Gamma$ finitely generated. Assume that $M$ is stably embedded in $\str G \M \p$. What can we say? 
\end{problem}

\noindent
Comment:
We know that when $\Gamma$ is finite then the extension can be interpreted in $M$ (Theorem \ref{ominimal}) and $G$ can be given a group topology making $\pi$ into a topological covering (Proposition \ref{real-case}). This is not always true when $\Gamma$ is infinite (even for $H=\R/\Z$) due to the non-categoricity result in Theorem \ref{modZ}. However one can at least conjecture that there is a topological covering homomorphism
$$\exact \Gamma {\i'} {G'} {\p'} H$$ such that $\str G \p M \equiv \str {G'} {\p'} M$.
This is a kind of converse to the statement that if $\exact \Gamma \i G \p H$ is a covering homomorphism, then $M$ is stably embedded in $\str G \M \p$ at the level of theories. The latter statement is implicit in earlier results. In fact we have already remarked that when $G$ is the universal cover Propositions \ref{coc} and \ref{uniform} apply, and the same argument works for every cover.

\end{document}